\documentclass[times]{elsarticle}

\usepackage{amsmath}
\usepackage{amsthm}
\usepackage{amssymb}
\usepackage{bm}
\journal{arXiv.org}

\begin{document}

\begin{frontmatter}

\title{Numerical investigation of a space-fractional model of turbulent fluid flow in rectangular ducts\tnoteref{label1}}
\tnotetext[label1]{This work was supported by the Russian Foundation for Basic Research (projects 14-01-00785, 15-01-00026).}

\author{Alexander~G.~Churbanov\fnref{lab1}}
\ead{chur@ibrae.ac.ru}

\author{Petr N. Vabishchevich\corref{cor1}\fnref{lab1,lab2}}
\ead{vabishchevich@gmail.com}
\cortext[cor1]{Correspondibg author.}

\address[lab1]{Nuclear Safety Institute, Russian Academy of Sciences,
              52, B. Tulskaya, 115191 Moscow, Russia}

\address[lab2]{North-Eastern Federal University,
	      58, Belinskogo, 677000 Yakutsk, Russia}

\begin{abstract}

The models that are based of fractional derivatives
should be highlighted among promising new models to describe turbulent fluid flows.
In the present work, a steady-state flow in a duct is considered under the condition
that the turbulent diffusion is governed by a fractional power of the Laplace operator.
To study numerically flows in rectangular channels, finite-difference approximations are employed.
For approximate solving the corresponding boundary value problem,
the iterative method of conjugate gradients is used. At each iteration,
the problem with a fractional power of the grid Laplace operator is solved.
Predictions of turbulent flows in ducts at different Reynolds numbers
are presented via mean velocity fields.

\end{abstract}

\begin{keyword}
Turbulent flow \sep fluid flow in ducts \sep fractional power of the Laplace operator \sep finite-difference problem \sep iterative method of conjugate gradients

\end{keyword}

\end{frontmatter}

\section{Introduction}

To model continuum mechanics phenomena, different models of turbulence are employed 
(see, e.g., \cite{book:9106,book:741002} among others). 
In terms of practical use, emphasis is on simple mathematical models of turbulence, 
which, on the one hand, are not much more complex in comparison with models for laminar flows, 
and on the other hand, reproduce the basic features of turbulent regimes of liquid and gas flows.

Nowadays, non-local applied mathematical models based on using fractional derivatives in time and space 
are actively discussed \cite{baleanu2012fractional,eringen2002nonlocal,kilbas2006theory}. 
Many models in applied physics, biology, hydrology and finance,
involve both sub-diffusion (fractional in time) and supper-diffusion (fractional in space) operators. 
Supper-diffusion problems are treated as evolutionary problems with a fractional power of an elliptic operator.

Such an anomalous diffusion model is used in \cite{book:1123818} to describe turbulent flows. 
In the work \cite{chen2006speculative}, a turbulent diffusion in the Reynolds equations for the mean velocity 
is governed by the fractional Laplacian. The development of this approach is
hindered by the lack of simple and robust numerical algorithms for solving 
boundary value problems for equations with fractional powers. 
In the best case (see, e.g., \cite{sousa2013}), investigations are restricted to simple one-dimensional
in space models.

For solving problems with fractional powers of elliptic operators, we can apply
finite volume and finite element methods oriented to using arbitrary
domains and irregular computational grids \cite{KnabnerAngermann2003,QuarteroniValli1994}.
The numerical implementation involves the matrix function-vector multiplication.
For such problems, different approaches \cite{higham2008functions} are available.
Problems of using Krylov subspace methods with the Lanczos approximation
when solving systems of linear equations associated with
the fractional elliptic equations are discussed in \cite{ilic2009numerical}.
A comparative analysis of the contour integral method, the extended Krylov subspace method, 
and the preassigned poles and interpolation nodes method for solving
space-fractional reaction-diffusion equations is presented in \cite{burrage2012efficient}.
The simplest variant is associated with the explicit construction of the solution using the known
eigenvalues and eigenfunctions of the elliptic operator with diagonalization of the corresponding matrix
\cite{bueno2012fourier,ilic2005numerical}. 
Unfortunately, all these approaches demonstrates too high computational complexity for multidimensional problems.

We have proposed \cite{vabishchevich2014numerical} a computational algorithm for solving
an equation with fractional powers of elliptic operators on the basis of
a transition to a pseudo-parabolic equation.
For the auxiliary Cauchy problem, the standard two-level schemes are applied.
The computational algorithm is simple for practical use, robust, and applicable to solving
a wide class of problems. A small number of time steps is required to find a solution.
This computational algorithm for solving equations with fractional powers of operators
is promising when considering transient problems. 

To implement numerically a space-fractional model of turbulent fluid flow,
we must take into account a multi-term structure of the problem operator.
Namely, here one term is the standard elliptic operator (normal diffusion), whereas the second term is
a fractional power of an elliptic operator (anomalous diffusion). 
For solving such non-classical problems, it seems natural to apply iterative methods 
with an appropriate choice of preconditioners \cite{book:1053511,book:8982}. 

In this paper, for predicting a steady-state turbulent flow in a duct,
we apply a model with a turbulent space-fractional diffusion. To solve
numerically this problem with the multi-term diffusion,
we employ the iterative method of conjugate gradients, where the problem with normal diffusion
is solved at each iteration to construct a preconditioner. For solving the problem with
the fractional Laplacian, a pseudo-parabolic equation is used.
The paper is organized as follows.
In Section 2, a mathematical model with the fractional Laplacian is introduced to describe a turbulent 
flow in a rectangular duct.
The discrete problem and computational algorithm are discussed in Section 3.
Section 4 presents an analysis of the impact of
the basic parameters of the problem on numerical results obtained using the developed model.

\section{A space-fractional model of turbulent fluid flow}

Motion of an incompressible fluid is governed by the Navier-Stokes equations:
\begin{equation}\label{1}
 \frac{\partial \bm v}{\partial t} + \bm v \cdot \nabla v + \frac{1}{\rho} \nabla p - 
 \nu \triangle \bm v = 0 ,
\end{equation} 
\begin{equation}\label{2}
 \nabla \cdot \bm v = 0 .
\end{equation}
Here $\rho$ is the density, $p$ denotes the pressure, $\bm v$ stands for the velocity vector,
and $\nu$ is the fluid viscosity. 

To obtain the Reynolds equations for turbulent flows \cite{book:9106,book:741002},
the velocity and pressure $\bm v, p$ are decomposed into the sum of the mean flow components
$\overline{\bm v}, \overline{p}$ and fluctuating components 
$\widetilde{\bm v}, \widetilde{p}$. Substituting this decomposition into (\ref{1}), (\ref{2}), 
we arrive at the  Reynolds equations written in the following coordinate-wise representation
($\bm v = (v_1,v_2, v_3)$):
\begin{equation}\label{3}
 \frac{\partial \overline{v}_i}{\partial t} + 
 \overline{v}_j \frac{\partial \overline{v}_i}{\partial x_j} + 
 \frac{1}{\rho} \frac{\partial \overline{p}}{\partial x_i} -
 \nu \triangle \overline{v}_i +
 \frac{\partial}{\partial x_i} \overline{\widetilde{v}_i \widetilde{v}_j} = 0 ,
\end{equation} 
\begin{equation}\label{4}
 \frac{\partial \overline{v}_i}{\partial x_i} = 0 .
\end{equation}
A RANS model of turbulence is defined by a particular formulation for the Reynolds stress tensor
${\displaystyle \rho \overline{\widetilde{v}_i \widetilde{v}_j}}$.
 
For the space-fractional model, we have
\begin{equation}\label{5}
 \frac{\partial}{\partial x_i} \overline{\widetilde{v}_i \widetilde{v}_j} = 
 \xi (-\triangle)^{\alpha} \overline{v}_i .
\end{equation}
Here the coefficient $\xi$ is treated as the eddy (turbulent) diffusivity. 
In the work \cite{chen2006speculative}, some  arguments are given in favor of setting 
the power $\alpha$ equal to $1/3$.
In view of (\ref{5}), equations (\ref{3}), (\ref{4}) may be written in the form similar to (\ref{1}), (\ref{2}),
i.e.,
\begin{equation}\label{6}
 \frac{\partial \overline{\bm v}}{\partial t} + 
 \overline{\bm v}\cdot \nabla \overline{\bm v} + \frac{1}{\rho} \nabla \overline{p} - 
 \nu \triangle \overline{\bm v}  +
 \xi (-\triangle)^{\alpha} \overline{\bm v} = 0 ,
\end{equation} 
\begin{equation}\label{7}
 \nabla \cdot \overline{\bm v} = 0 .
\end{equation} 

Let us consider a steady-state stabilized in the longitudinal direction flow in rectangular channels
($\bm x = (x_1,x_2)$):
\[
 \Omega = \{ \bm x  \ | \ \bm x = (x_1,x_2) \ 0 < x_i < d_i, \ i = 1,2 \} .
\]
Let $x_3$ be the longitudinal coordinate and assume that $\overline{\bm v} = (0,0,u)$. 
Then from (\ref{6}), (\ref{7}), we obtain the following equation for the longitudinal component of the velocity:
\begin{equation}\label{8}
 - \nu \triangle u +
 \xi (-\triangle)^{\alpha} u = \chi  ,
 \quad  \bm x \in \Omega ,
\end{equation}  
where
\[
 p = p(x_3),
 \quad  \chi  = - \frac{1}{\rho} \frac{d p}{d x_3}.  
\]
The equation (\ref{8}) is supplemented with homogeneous Dirichlet boundary conditions:
\begin{equation}\label{9}
 u(\bm x) = 0, 
 \quad \bm x  \in \partial \Omega ,
\end{equation}
which corresponds to the no-slip condition on rigid walls.

For the normalization of equation (\ref{8}), as the reference values, we employ the channel height 
$d_2$ and the velocity scale
\[
 u_0 =  \frac{d_2^2}{\nu } \, \chi .
\]
For the dimensionless velocity $u$, using for the dimensionless quantities 
the same notation as for the dimensional ones, we obtain
\begin{equation}\label{10}
 - \triangle u +
 \mu  (-\triangle)^{\alpha} u = 1  ,
 \quad  \bm x \in \Omega ,
\end{equation}  
where 
\[
 \Omega = \{ \bm x \ | \ \bm x = (x_1,x_2) \ | \ 0 < x_1 < d, \  0 < x_2 < 1 \} ,
\]
\[
 \mu = \frac{\xi }{\nu } \, d_2^{2(1-\alpha)} . 
\] 
Thus, the boundary value problem (\ref{9}), (\ref{10}) has three governing parameters, namely, $\alpha, \mu$ and $d$. 

\section{Computational algorithm} 

To solve the steady-state problem  (\ref{9}), (\ref{10}), we introduce a uniform grid in the domain $\Omega$:
\[
\overline{\omega}  = \{ \bm{x} \ | \ \bm{x} =\left(x_1, x_2\right), \quad x_k =
i_k h_k, \quad i_k = 0,1,...,N_k,
\quad N_1 h_1 = d,  \quad N_2 h_2 = 1 \} ,
\]
with $\overline{\omega} = \omega \cup \partial \omega$, where
$\omega$ is the set of interior points and $\partial \omega$ is the set of boundary grid points. 
For grid functions $y(\bm x)$ such that $y(\bm x) = 0, \ \bm x \notin \omega$, we define the Hilbert space 
$H=L_2\left(\omega\right)$, where the scalar product and the norm are given as follows:
\[
\left(y, w\right) \equiv  \sum_{\bm x \in  \omega} y\left(\bm{x}\right)
w\left(\bm{x}\right) h_1 h_2,  \quad 
\| y \| \equiv  \left(y, y\right)^{1/2}.
\]

For the discrete Laplace operator $A$, we introduce the additive representation
\begin{equation} \label{11}
A = \sum_{k=1}^{2} A_k, 
\quad \bm{x} \in \omega,
\end{equation}
where $A_k, \ k=1,2$ are  associated with the corresponding differential operator of the second derivative 
in one direction.

For all grid points except adjacent to the boundary, the grid operator $A_1$ 
can be written as
\[
  \begin{split}
  A_1 y = & -
  \frac{1}{h_1^2} (y(x_1+h_1,h_2) - 2 y(\bm{x}) + y(x_1-h_1,h_2)), \\
  & \qquad \bm{x} \in \omega, 
  \quad x_1 \neq 0.5h_1,
  \quad x_1 \neq d-0.5h_1. 
 \end{split} 
\]
In the points that are adjacent to the boundary, the approximation is constructed taking into account 
the boundary condition (\ref{9}):
\[
  \begin{split}
  A_1 y & =  -
  \frac{1}{h_1^2} (y(x_1+h_1,h_2) - 2 y(\bm{x})) , \quad \bm{x} \in \omega, 
  \quad x_1 = 0.5h_1, \\
  A_1 y & = 
  \frac{1}{h_1^2} (2 y(\bm{x}) - y(x_1-h_1,h_2)), 
  \quad \bm{x} \in \omega, 
  \quad x_1 = d-0.5h_1. 
 \end{split} 
\]
Similarly we construct the grid operator $A_2$.
For the above grid operators, we have (see, e.g., \cite{Samarskii,SamarskiiNikolaev1978})
\[
  A_k = A_k^* \geq \delta_k E, 
  \quad  \delta_k = \frac{4}{h_k^2} \sin^2 \frac{\pi }{2 N_k} , 
  \quad  k=1,2 . 
\]
where $E$ is the identity operator.
Because of this, the discrete Laplace operator (\ref{11}) is self-adjoint and positive definite in $H$:
\begin{equation} \label{12}
 A = A^* \geq \delta E,
 \quad \delta =  \sum_{k=1}^{2} \delta_k.
\end{equation}
It approximates the differential Laplace operator with the truncation error
$\mathcal{O} \left(|h|^2\right)$, $|h|^2 = h_1^2+h_2^2$. 

To handle the fractional power of the grid operator $A$, let us consider the eigenvalue problem
\[
 A \varphi_m = \lambda_m \varphi_m , 
\]
which has the well-known analytical solution. We have
\[
 \delta = \lambda_1 \leq \lambda_2   \leq ... \leq \lambda_M,
 \quad M = (N_1-1)(N_2-1) , 
\]
where eigenfunctions $\varphi_m, \ \|\varphi_m\| = 1, \ m = 1,2, ..., M$
form a basis in $H$. Therefore
\begin{equation}\label{13}
 y = \sum_{m= 1}^{M}(y, \varphi_m) \varphi_m . 
\end{equation} 
For the fractional power of the operator $A$, we have
\[
 A^\alpha y = \sum_{m= 1}^{M}(y, \varphi_m) \lambda_m^\alpha \varphi_m .
\] 

Using the above approximations, we arrive from (\ref{9}), (\ref{10}) at the discrete problem
\begin{equation}\label{14}
  A y + \mu A^\alpha y = 1 .
\end{equation}
In our particular case with using uniform meshes in a rectangle, the solution of equation (\ref{14}) 
can be constructed explicitly via the known eigenvalues and eigenfunctions. For the solution represented
in the form of (\ref{13}), we obtain
\[
 (y, \varphi_m) = \frac{(1, \varphi_m)}{\lambda_m + \mu \lambda_m^\alpha} ,
 \quad   m = 1,2, ..., M.
\]
We are interested in solving problems of type (\ref{14}) under more general conditions, 
where the complete eigenvalue problem requires large computational costs.

In this situation, we cannot directly apply the well-developed iterative methods of linear algebra
and an appropriate software for solving (\ref{14}). This results from two reasons.
On the one hand, we have the term $\mu A^\alpha y$ on the left side. On the other hand, the equation is multi-term, 
i.e., it is represented as the sum of two individual operators.

Obviously, elliptic problems with the operator $A$ can be solved in an efficient way.
Then the operator $A$ can be selected as a preconditioner for iterative solving equation (\ref{14}).
Let $y_k$ be an approximate solution at the $k$-th iteration. 
If we apply the conjugate gradient method \cite{book:1053511,book:8982},
then the new iteration is defined as follows. Denote 
$r_k = 1 - \widetilde{A} y_k$, $\widetilde{A} = A + \mu A^\alpha$  as the original residual and 
let $z_k = A^{-1} r_k$ be the residual for the preconditioned equation. With the initial $p_0 = z_0$
and the given $y_0$, for $k = 0, 1, ...$, we have
\begin{equation}\label{15}
\begin{split}
 \alpha_k & = \frac{(z_k, r_k)}{(\widetilde{A} p_k, p_k)} ,
 \quad y_{k+1} = y_k + \alpha_k p_k,
 \quad r_{k+1} = r_k - \alpha_k \widetilde{A} p_k,  \\
 z_{k+1} & = A^{-1} r_{k+1},
 \quad  \beta_k = \frac{(z_{k+1}, r_{k+1})}{(z_k, r_k)} ,
 \quad p_{k+1} = z_{k+1} + \beta_k p_k . 
\end{split}
\end{equation} 

The convergence rate of the iterative method (\ref{15})
is governed \cite{SamarskiiNikolaev1978} by the constants $\gamma_1$ and $\gamma_2$
(more precisely, by the ratio $\varkappa = \gamma_1 / \gamma_2$) 
in the following bilateral operator inequality:
\begin{equation}\label{16}
 \gamma_1 A \leq A + \mu A^\alpha \leq \gamma_2 A ,
 \quad \gamma_1 > 0 .
\end{equation} 
For $\widetilde{A}$,  in view of (\ref{12}) and $0 < \alpha < 1$, we have
\[
 \widetilde{A} = A + \mu A^\alpha > A,
 \quad A + \mu A^\alpha = (E + \mu A^{\alpha-1} ) A \leq 
 (1 + \mu \delta^{\alpha-1}) A . 
\] 
Therefore for $\gamma_1$ and $\gamma_2$ in (\ref{16}), we obtain
\[
  \gamma_1 = 1,
  \quad \gamma_2 =  1 + \mu \delta^{\alpha-1} .
\]
This establishes the dependence of the number of iterations in the conjugate gradient method (\ref{16})
on $\mu, \ \delta$ and $\alpha$.

At each iteration, we must evaluate the quantity
\[
 \widetilde{A} p_k = A (p_k + \mu A^{\alpha-1} p_k) .
\]
The emphasis here is on calculating $w = A^{\alpha-1} p_k$. 
It is necessary to solve the problem
\begin{equation}\label{17}
 A^\beta w = f,
\end{equation} 
where $\beta = 1-\alpha, \ f = p_k$ for $0 < \beta < 1$.
We apply the approach proposed in the paper \cite{vabishchevich2014numerical}.

An approximate solution is sought as the solution of an auxiliary evolutionary problem, 
where $t$ is the pseudo-time evolution variable. Assume that
\[
 v(t) = (\theta \delta)^{\alpha} (t (A - \theta \delta E) + \theta  \delta E)^{-\alpha} v(0) ,
\]
with $0 < \theta < 1$. Therefore
\[
 v(1) =  (\theta \delta)^{\alpha} A^{-\alpha} v(0) 
\]  
and then $w = v(1)$.
The function $v(t)$ satisfies the evolutionary equation
\begin{equation}\label{18}
  (t D + \theta \delta I) \frac{d v}{d t} + \alpha D v = 0 ,
  \quad 0 < t \leq 1 ,
\end{equation}  
where
\[
 D = A - \theta \delta E .
\] 
By (\ref{12}), we get
\begin{equation}\label{19}
 D = D^* \geq (1-\theta) \delta E > 0 .
\end{equation}
We supplement equation (\ref{18}) with the initial condition
\begin{equation}\label{20}
 v(0) = (\theta \delta)^{-\alpha} f .  
\end{equation} 
The solution of equation (\ref{17}) can be defined as the solution
of the Cauchy problem (\ref{18})--(\ref{20}) at the final time moment $t=1$.
In \cite{vabishchevich2014numerical}, the case of $\theta =1$ was studied.

For the solution of the problem (\ref{18}), (\ref{20}), we can obtain various a priori estimates.
Elementary estimates have the form
\begin{equation}\label{21}
  \|v(t)\|_G \leq \|v(0)\|_G , 
\end{equation} 
where, for instance, $G = E, D$.
To obtain (\ref{21}) for $G = D$, multiply scalarly equation (\ref{18}) by $dv/ dt$.
For $G = E$,  equation (\ref{18}) is multiplied by $\alpha v + t dv/ dt$.

To solve numerically the problem (\ref{18}), (\ref{20}),
we use a simple two-level scheme.
Let $\tau$ be a step of a uniform grid in time such that $v^n = v(t^n), \ t^n = n \tau$,
$n = 0,1, ..., N_0, \ N_0\tau = 1$.
Let us approximate equation (\ref{18}) by the Crank-Nicolson scheme
\begin{equation}\label{22}
 (t^{n+1/2} D + \theta \delta E) \frac{ v^{n+1} - v^{n}}{\tau }
 + \alpha D  \frac{ v^{n+1} + v^{n}}{2} = 0,  \quad n = 0,1, ..., N_0-1,
\end{equation}
\begin{equation}\label{23}
 v^0 = (\theta \delta)^{-\alpha} f .
\end{equation} 
The difference scheme (\ref{22}), (\ref{23}) approximates the problem 
(\ref{18}), (\ref{20})
with the second order by $\tau$.

For $0 < \theta < 1$,  the difference scheme (\ref{22}), (\ref{23}) is unconditionally stable
with respect to the initial data.
The approximate solution satisfies the estimate 
\begin{equation}\label{24}
  \|v^{n+1}\|_G \leq \|v^0\|_G , 
  \quad n = 0,1, ..., N_0-1,
\end{equation} 
with $G = E, D$.

Multiplying scalarly equation (\ref{22}) by $v^{n+1} - v^{n}$, we get
\[
  \|v^{n+1}\|_D \leq \|v^n\|_D , 
  \quad n = 0,1, ..., N_0-1 .  
\]
This inequality ensures the estimate (\ref{24}) for $G = D$.

Similarly, we consider the case with $G = E$.
Rewrite equation (\ref{22}) in the form
\[
 \theta \delta \frac{ v^{n+1} - v^{n}}{\tau } + 
 D \left( \alpha \frac{ v^{n+1} + v^{n}}{2}+ t^{n+1/2}\frac{ v^{n+1} - v^{n}}{\tau } \right ) = 0 . 
\] 
Multiplying scalarly it by
\[
 \alpha \frac{ v^{n+1} + v^{n}}{2}+ t^{n+1/2}\frac{ v^{n+1} - v^{n}}{\tau },
\] 
in view of (\ref{19}), we arrive at
\[
 \left ( \frac{ v^{n+1} - v^{n}}{\tau }, \frac{ v^{n+1} + v^{n}}{2} \right ) \leq 0 . 
\] 
We have
\[
  \|v^{n+1}\| \leq \|v^n\| , 
  \quad n = 0,1, ..., N_0-1 , 
\] 
that means the fulfilment of the estimate (\ref{24}) with $G = E$.

\section{Numerical results} 

\begin{figure}
  \begin{center}
    \includegraphics[width=0.7\linewidth] {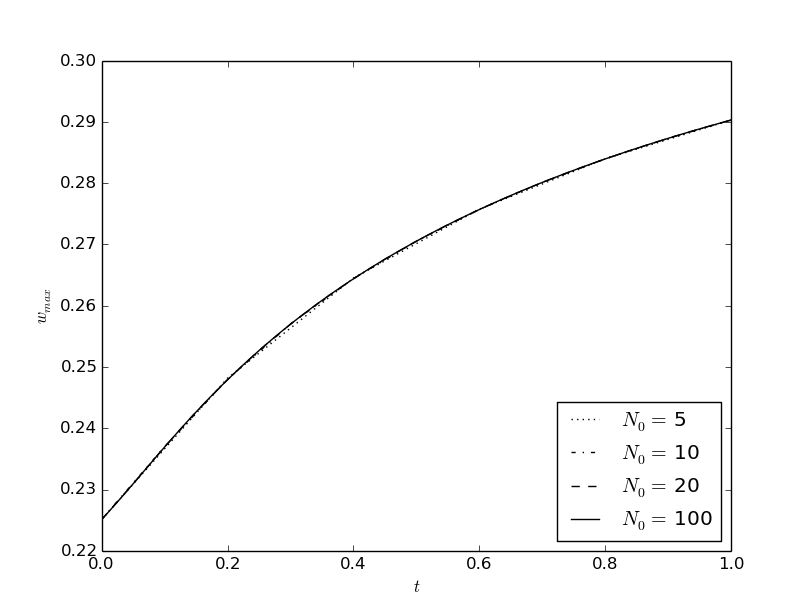}
	\caption{Evolution histories of $w_{max}$ for different grids in $t$ ($\beta = 0.5$, $\theta \delta = 2 \pi^2$)}
	\label{f-1}
  \end{center}
\end{figure}

\begin{figure}
  \begin{center}
    \includegraphics[width=0.7\linewidth] {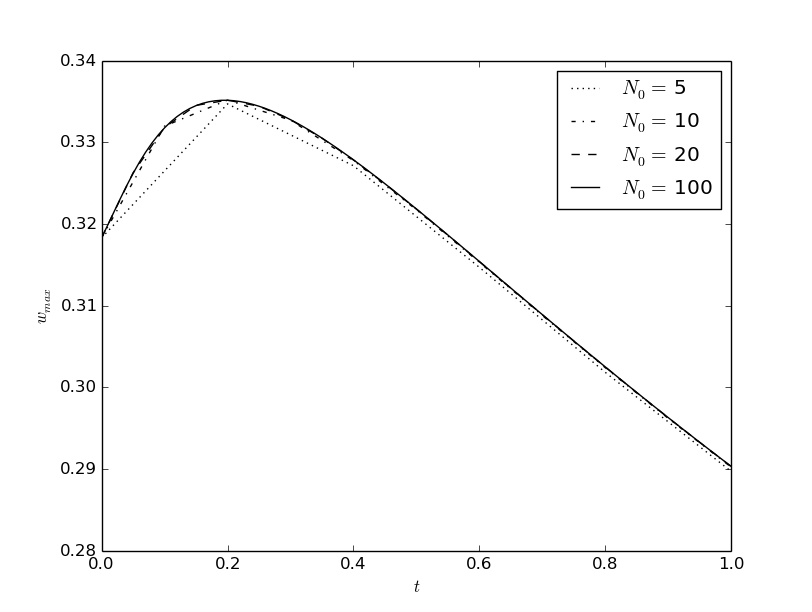}
	\caption{Evolution histories of $w_{max}$ for different grids in $t$ ($\beta = 0.5$, $\theta \delta = \pi^2$)}
	\label{f-2}
  \end{center}
\end{figure}

To discuss our predictions, we start with calculations of the problem with the fractional power (\ref{17}).
The scheme (\ref{22}), (\ref{23}) was applied. The problem, unless otherwise stated, was solved
on the spatial grid $N_1 = N_2 = 100$, $d = 1$ with $\beta = 0.5$, $\theta \delta = 2 \pi^2$.
The evolution histories of the maximum value $w_{max}$ of the approximate solution 
(located at the center of the domain) are shown in Figure~\ref{f-1} for 
various computational grids in the pseudo-time evolution variable $t$ ($N_0 = 5, 10, 20, 100$).
It is easy to see that even on coarse grids in $t$, we observe a good accuracy of the solution.
Figure~\ref{f-2} demonstrates similar pseudo-time histories calculated starting from other initial value of (\ref{23})
in the scheme (\ref{22}). Namely, in this case, we use a rougher initial approximation $\theta \delta = \pi^2$
that corresponds to an inaccurate estimation for the lower bound (\ref{12}) of the operator $A$.

\begin{figure}
  \begin{center}
    \includegraphics[width=0.7\linewidth] {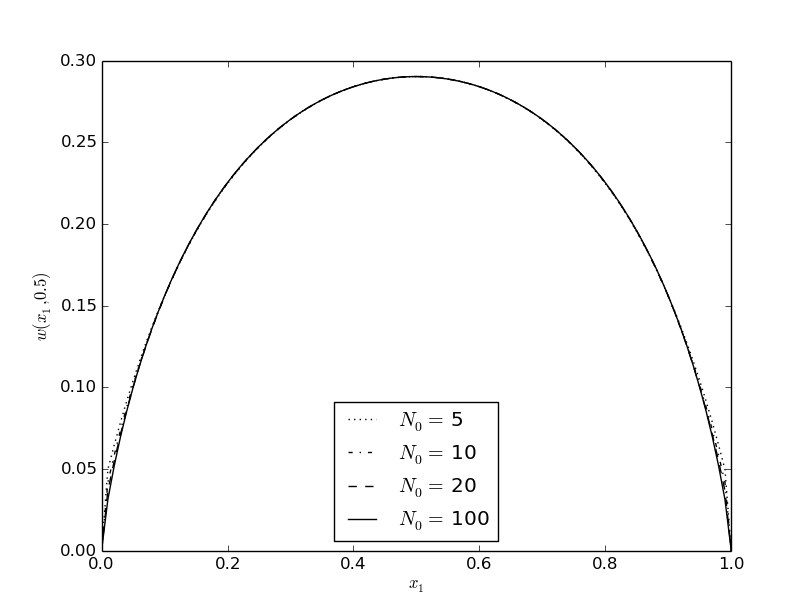}
	\caption{Profiles of $w(x_1,0.5)$ for different grids in $t$ ($\beta = 0.5$, $\theta \delta = 2 \pi^2$)}
	\label{f-3}
  \end{center}
\end{figure}

\begin{figure}
  \begin{center}
    \includegraphics[width=0.7\linewidth] {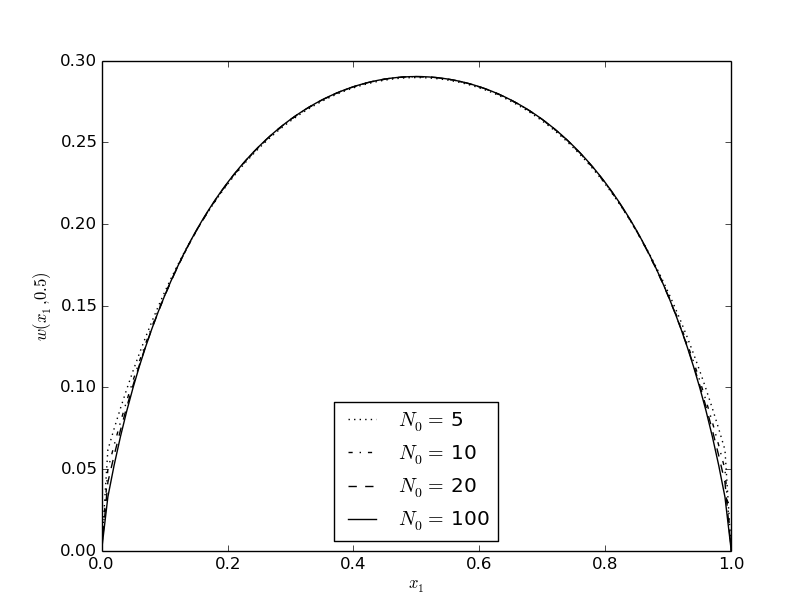}
	\caption{Profiles of $w(x_1,0.5)$ for different grids in $t$ ($\beta = 0.5$, $\theta \delta = \pi^2$)}
	\label{f-4}
  \end{center}
\end{figure}

The non-local convergence of the approximate solution with refining the grid in the pseudo-time evolution variable $t$
is depicted in Figures~\ref{f-3} and \ref{f-4} for solution profiles. There are presented profiles
of the solution of the problem (\ref{12}) in the mid-section ($x_2 = 0.5$). As above, these profiles of $w(x_1,0.5)$ were
calculated using various grids in $t$ and starting with two different initial values. 
Obviously, the maximum error is observed in the vicinity of boundaries of the computational domain.



\begin{figure}
  \begin{center}
    \includegraphics[width=0.7\linewidth] {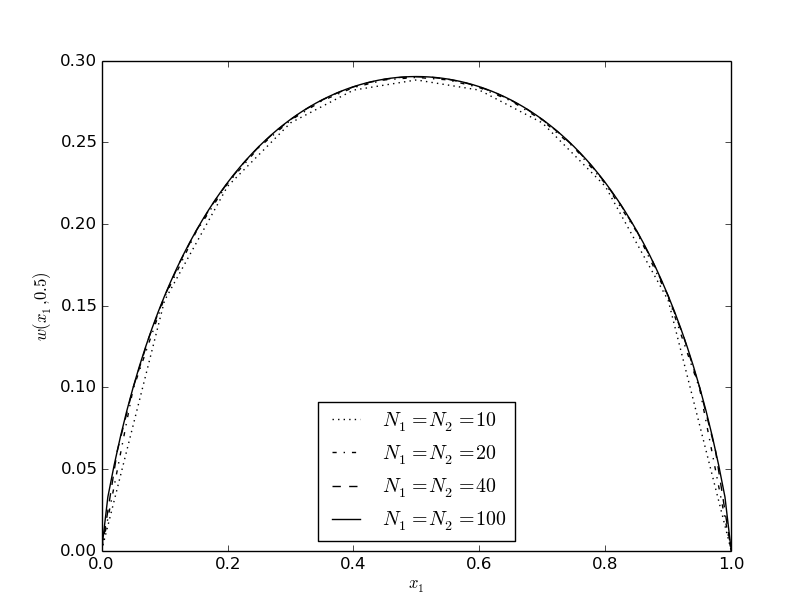}
	\caption{Profiles of $w(x_1,0.5)$ for different grids in space ($\beta = 0.5$)}
	\label{f-5}
  \end{center}
\end{figure}

The solution convergence for the fractional Laplace operator problem  with refining the grid in space is shown 
in Figure~\ref{f-5} for the above mid-section profiles of $w(x_1,0.5)$.
The calculations were performed on the finest grid in $t$ ($N_0 = 100$).

\begin{figure}
  \begin{center}
    \includegraphics[width=0.7\linewidth] {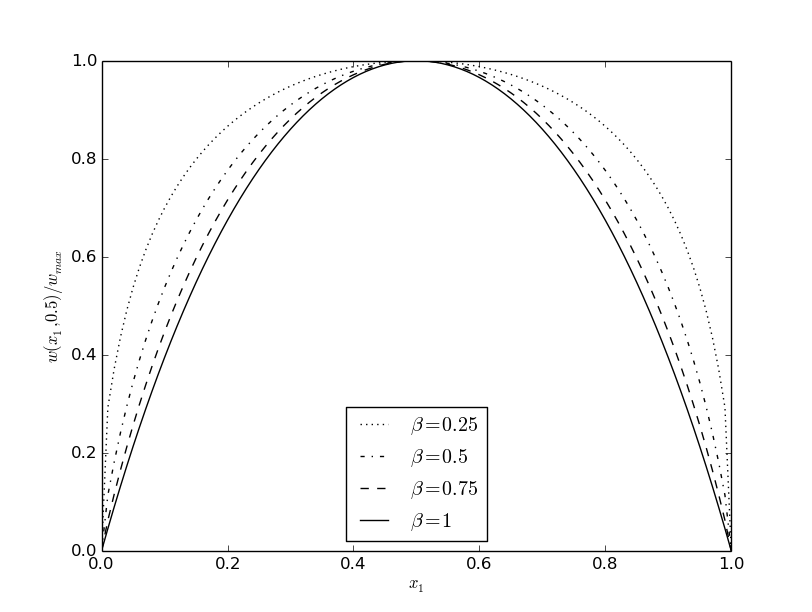}
	\caption{Profiles of $w(x_1,0.5)$ for various values of the power $\beta$}
	\label{f-6}
  \end{center}
\end{figure}

\begin{figure}
  \begin{center}
    \includegraphics[width=0.7\linewidth] {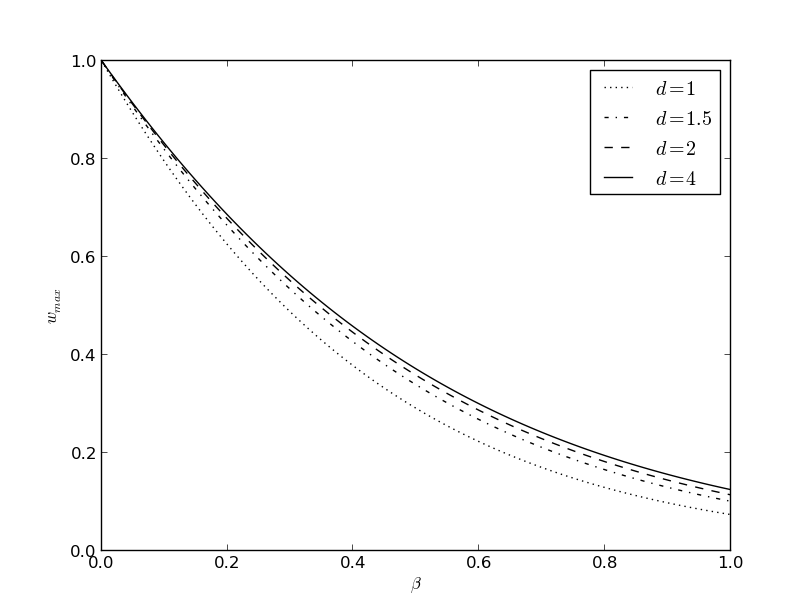}
	\caption{The solution maximum $w_{max}$ for various values of power $\beta$ and width $d$}
	\label{f-7}
  \end{center}
\end{figure}

For the problem (\ref{17}), the main interest is in the impact of the power $\beta$ on the solution features. 
To eliminate the influence of grid parameters (grid steps in space and the pseudo-time evolution variable $t$),
all predictions in this parametric study were performed on the finest grid with $N_1 = N_2 = 100$ and $N_0 = 100$.
Figure~\ref{f-6} presents mid-section profiles of the solution for various values of $\beta$.
For the convenience of a comparison, the solutions are normalized to the maximum value. 
It is easy to see that the decreasing of $\beta$ leads to more gently sloping profiles. 
When $\beta \rightarrow 0$, we have $w \rightarrow f$ in the computational domain $\Omega$. 
The dependence of the solution maximum $w_{max}$ on the power $\beta$ and the geometry 
(the width of the computational domain $d$) is presented in Figure~\ref{f-7}.

\begin{figure}
  \begin{center}
    \includegraphics[width=0.7\linewidth] {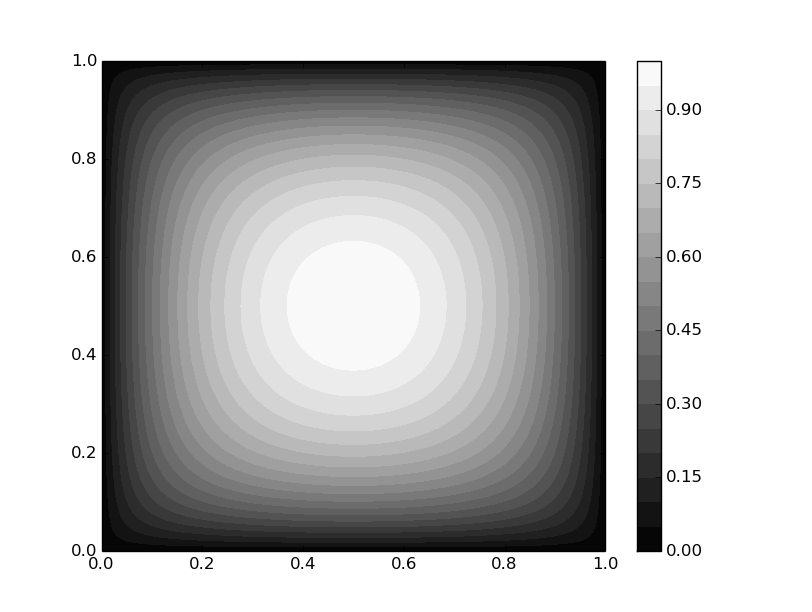}
	\caption{Isocontoures of the solution for $\beta = 0.75$}
	\label{f-8}
  \end{center}
\end{figure}

\begin{figure}
  \begin{center}
    \includegraphics[width=0.7\linewidth] {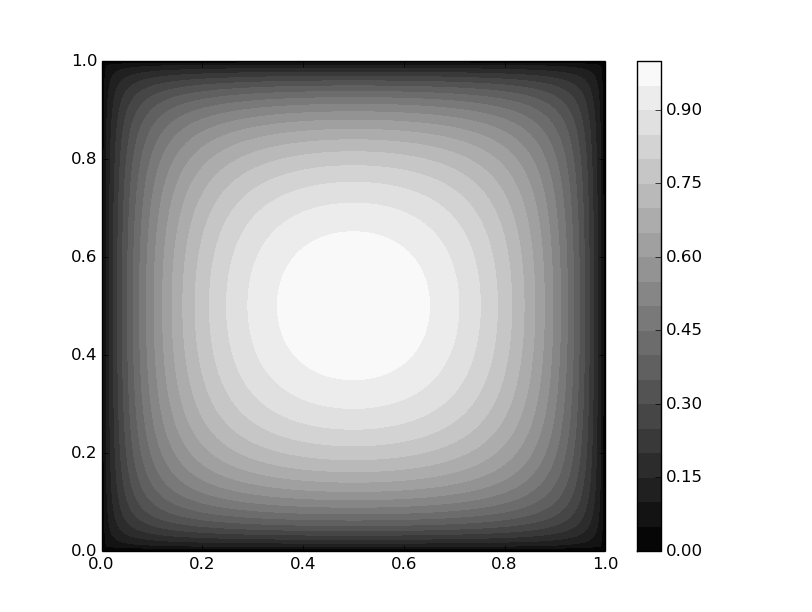}
	\caption{Isocontoures of the solution for $\beta = 0.5$}
	\label{f-9}
  \end{center}
\end{figure}

\begin{figure}
  \begin{center}
    \includegraphics[width=0.7\linewidth] {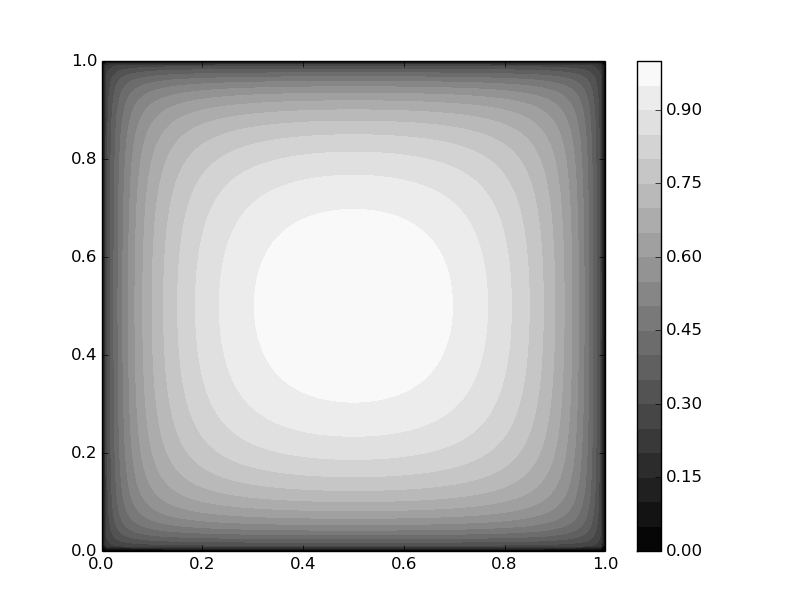}
	\caption{Isocontoures of the solution for $\beta = 0.25$}
	\label{f-10}
  \end{center}
\end{figure}

The solution of the problem (\ref{17}) normalized to the maximum value is shown in  Figures~\ref{f-8}--\ref{f-10}
as isocontoures in the whole computational domain for different values of the power $\beta$. 
We can observe the formation of a boundary layer when $\beta \rightarrow 0$.

\begin{figure}
  \begin{center}
    \includegraphics[width=0.7\linewidth] {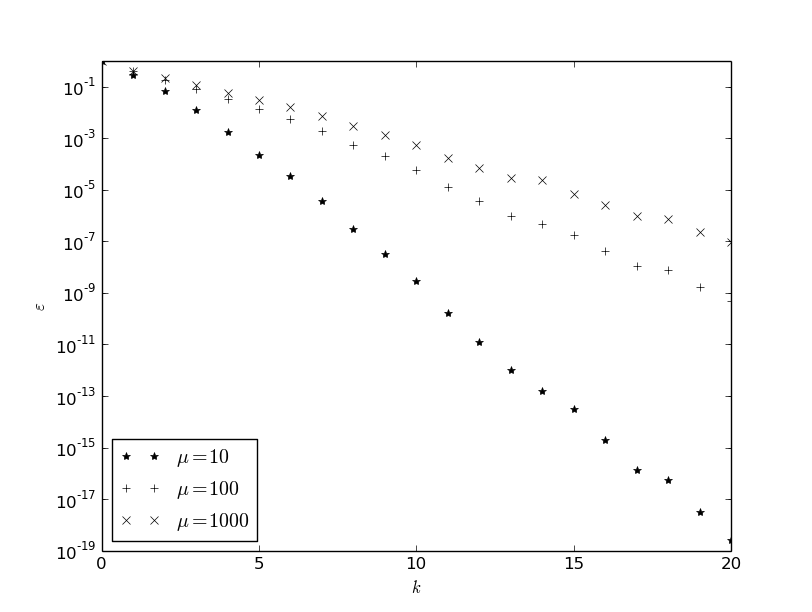}
	\caption{The relative error $\varepsilon$ vs iteration number $k$ for various values of $\mu$ ($\alpha = 0.5$)}
	\label{f-11}
  \end{center}
\end{figure}

\begin{figure}
  \begin{center}
    \includegraphics[width=0.7\linewidth] {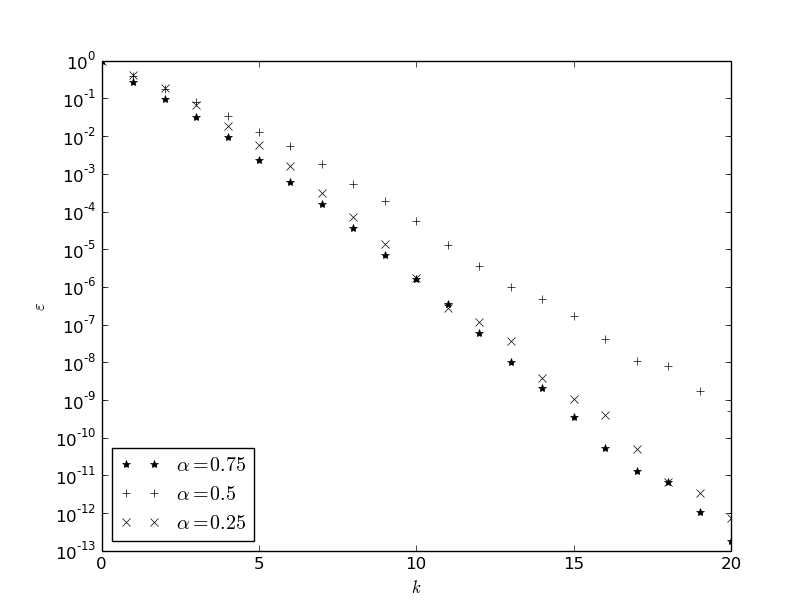}
	\caption{The relative error $\varepsilon$ vs iteration number $k$ for various values of $\alpha$ ($\mu = 100$)}
	\label{f-12}
  \end{center}
\end{figure}

Now we discuss the main object of our study, i.e., the problem (\ref{14}). To solve it, we apply 
the iterative method of conjugate gradients (\ref{15})  with the operator $A$ as a preconditioner. 
From the methodological point of view, the most interesting is the dependence of the iteration number
on the parameters $\mu$ and $\alpha$. The decreasing of
the relative error  $\varepsilon_k = \|r_k\| / \|r_0\|$ during iterations
(with the initial approximation $y_0 = 0$) is given in Figure~\ref{f-11}
for various values of $\mu$. The problem was solved with $\alpha = 0.5$. 
The dependence of the convergence rate on $\alpha$ for $\mu = 100$ is presented in Figure~\ref{f-12}.

\begin{figure}
  \begin{center}
    \includegraphics[width=0.7\linewidth] {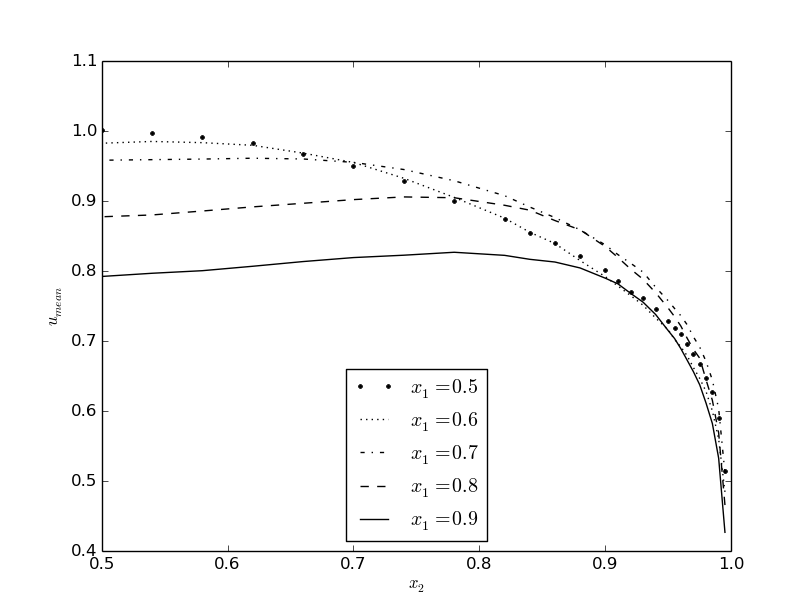}
	\caption{Experimental profiles of $u_{mean}$ \cite{1992forced} along various cross-lines of $x_1$}
	\label{f-13}
  \end{center}
\end{figure}

In modeling turbulent flows by means of the space-fractional model, we operate only with
mean values of the longitudinal velocity component. A more detailed description of turbulent flows
is carried out on the basis of more complicated models of turbulence (see, e.g., \cite{myong1991numerical,komen2014quasi}).
To validate our space-fractional model of turbulence, a comparison with experiments was done.
A fully developed turbulent flow in a square duct was measured in \cite{1992forced}. 
We use experimental data from this study, which are placed
on the Internet resource http://www.jsme.or.jp/ted/HTDB/fw.html.  
Experimental profiles of the normalized mean longitudinal velocity
$u_{mean}$ are shown in Figure~\ref{f-13} for various cross-sections of $x_1$ for a half of a cross-section.
Here the origin of coordinates is located at the left bottom
corner of the duct cross-section and so, at the center of the duct we have $x_1 = 0.5, x_2 = 0.5$. 
We see more gently sloping profiles of the velocity in approaching to duct walls. 
Also we see increasing of the velocity towards the corners of the duct, which is associated with secondary flows
observed in experiments and which it is difficult to reproduce using simple models of turbulence.

\begin{figure}
  \begin{center}
    \includegraphics[width=0.7\linewidth] {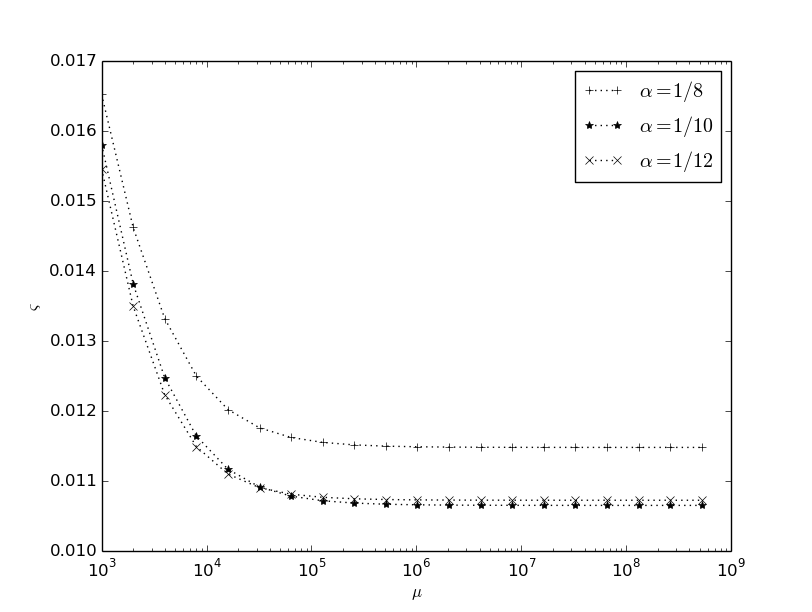}
	\caption{Dependence of $\varsigma$ on $\mu$ for optimal values of  $\alpha$}
	\label{f-14}
  \end{center}
\end{figure}

These experimental data we used to tune the parameters of our space-fractional model of turbulence (\ref{14})
in order to meet the above experimental data in the best way. For this purpose, a parametric study  
with respect to $\mu$ and $\alpha$ was done. We estimated the deviation between the calculated and measured values 
of the longitudinal velocity. Let $\bm x_l, \ l = 1,2, ..., L$ be the points of measurement. 
The deviation measure is the quantity
\[
 \varsigma(\mu, \alpha) = \frac{1}{L} \left (\sum_{l=1}^{L} (y(\mu, \alpha;\bm x_l) - u(\bm x_l) )^2 \right )^{1/2}, 
\]
where $y(\mu, \alpha;\bm x_l)$ is the predicted velocity, whereas $u(\bm x_l)$ is the measured velocity.
Figure~\ref{f-14} demonstrates the dependence of $\varsigma$ on $\mu$ for optimal values of $\alpha$. 
These results show that the first term in the left-hand side of equation (\ref{14}) can be neglected. 
Therefore it is  possible to use the one-term diffusion model, where instead of (\ref{6}), we consider the equation
\[
 \frac{\partial \overline{\bm v}}{\partial t} + 
 \overline{\bm v}\cdot \nabla \overline{\bm v} + \frac{1}{\rho} \nabla \overline{p} +
 \xi (-\triangle)^{\alpha} \overline{\bm v} = 0 .
\]
For a flow in a duct, we can reduce equation (\ref{8}) to the following equation
\begin{equation}\label{25}
 \xi (-\triangle)^{\alpha} u = \chi  ,
 \quad  \bm x \in \Omega ,
\end{equation} 
for the longitudinal velocity.

\begin{figure}
  \begin{center}
    \includegraphics[width=0.7\linewidth] {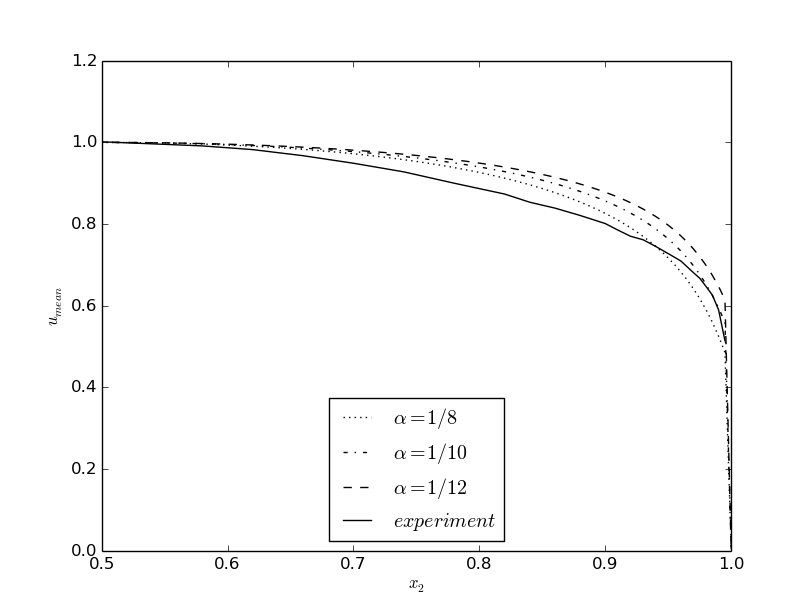}
	\caption{Comparison of experimental (solid line) and numerical profiles of $u_{mean}$ for $x_1 = 0.5$ }
	\label{f-15}
  \end{center}
\end{figure}

Numerical results obtained using the one-term space-fractional model of turbulence (\ref{25})
at near optimal values of $\alpha$ are presented in Figures~\ref{f-15}--\ref{17}.
The calculated data are compared with experimental profiles along three different lines of the duct
$x_1 = 0.5$, $x_1 = 0.7$ and $x_1 = 0.9$, respectively.
A good agreement between approximate solutions and the measurements is observed in the critical region near duct walls.
Relatively large discrepancies  take place only in the central zone of the duct cross-section.
This is partly due to the fact that the measurement points have a non-uniform distribution, i.e.,
near the boundaries the distance between the points is eight times lower than near the center.

\begin{figure}
  \begin{center}
    \includegraphics[width=0.7\linewidth] {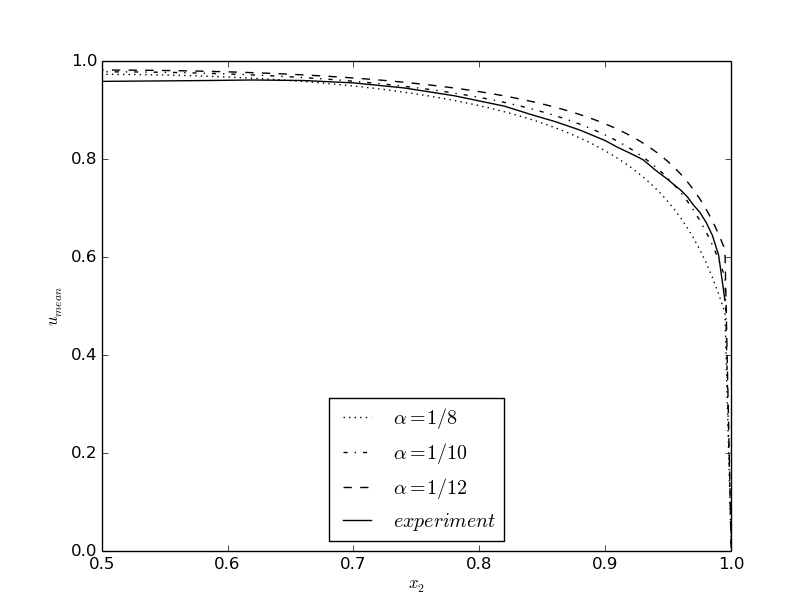}
	\caption{Comparison of experimental (solid line) and numerical profiles of $u_{mean}$ for $x_1 = 0.7$ }
	\label{f-16}
  \end{center}
\end{figure}

\clearpage

\begin{figure}
  \begin{center}
    \includegraphics[width=0.7\linewidth] {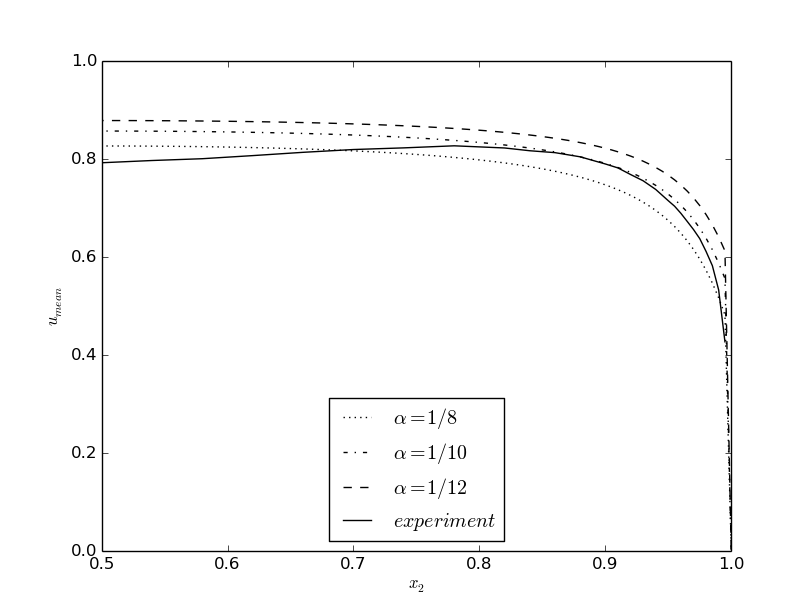}
	\caption{Comparison of experimental (solid line) and numerical profiles of $u_{mean}$ for $x_1 = 0.9$ }
	\label{f-17}
  \end{center}
\end{figure}

\end{document}